\documentclass[11pt,a4paper]{article}
\usepackage{amsmath, amssymb, theorem, latexsym,color}
\usepackage{graphicx}

\newtheorem{theorem}{Theorem}[section]
\newtheorem{lemma}[theorem]{Lemma}
\newtheorem{proposition}[theorem]{Proposition}
\newtheorem{definition}[theorem]{Definition}

\newtheorem{conjecture}[theorem]{Conjecture}

\newtheorem{remark}[theorem]{Remark}
\newtheorem{remarks}[theorem]{Remarks}

\newtheorem{property}[theorem]{Property}

\allowdisplaybreaks

\def\neweq#1{\begin{equation}\label{#1}}
\def\endeq{\end{equation}}

\def\Om {{\Omega}}
\def\la {{\lambda}}
\def\ep {{\epsilon}}

\def \qz { {\mathbb Q}}
\def \rz { {\mathbb R}}
 
\def\Sg {{\cal S}} 
\def \Ab{{\bf A}}
\def \Bb{{\bf B}}
\def \Xb{{\bf X}}
\def\qz {{\mathbb Q}}

\def \rz {{\mathbb R}}

\def \Inte{{\rm Int\,}}

\newcommand {\ar}{\rightarrow}
\newcommand {\pa}{\partial}

\numberwithin{equation}{section}

\title{A review on large $k$ minimal spectral  $k$-partitions and  Pleijel's Theorem.}

\author{ B. Helffer (Universit\'e Paris-Sud 11) \\
 and\\
 and T. Hoffmann-Ostenhof (University of Vienna)}
 \date{}

\begin{document}
\maketitle

 \begin{abstract}
In this survey, we review the properties
of minimal spectral $k$-partitions  in the two-dimensional case and revisit their connections with Pleijel's Theorem. We focus on the large $k$ problem (and the hexagonal conjecture) in connection
with two recent preprints by  J. Bourgain and S. Steinerberger on the Pleijel Theorem. This leads us also to discuss  some  conjecture by I. Polterovich, in relation with square tilings. We also establish a  Pleijel Theorem for Aharonov-Bohm Hamiltonians and  deduce from it, via the magnetic characterization of the minimal partitions, some lower bound for the number of critical points of a minimal partition.  \end{abstract}
\section{Introduction}

We consider  the Dirichlet  Laplacian in a bounded
regular domain $\Omega\subset \mathbb R^2$.
In \cite{HHOT} we have started to analyze the  relations between the nodal domains
 of the real-valued  eigenfunctions of this Laplacian and the partitions of
 $\Omega$ 
 by $k$
 open sets $ D_i$ which are minimal in the sense that  the
 maximum over the $ 
D_i$'s of the ground state energy\footnote{The ground state energy is
  the smallest eigenvalue.}  of  the Dirichlet
 realization
 of the Laplacian  in $ D_i$ is minimal. 
  We denote by $ \lambda_j(\Omega)$ the
increasing sequence of its eigenvalues and by $\phi_j$ some associated
 orthonormal basis of real-valued eigenfunctions. 
The groundstate  $ \phi_1$ can be chosen to be
strictly positive in $ \Om$, but the other eigenfunctions 
$  \phi_k$ must have zerosets.
For any real-valued continuous function $u$ on 
$ \overline\Om $, we define the zero set as 
\begin{equation}
N(u)=\overline{\{x\in \Om\:\big|\: u(x)=0\}}
\end{equation}
and call the components of $ \Om\setminus N(u)$ the nodal
domains of $u$.
The number of 
nodal domains of $ u$ is called $ \mu(u)$. These
$\mu(u)$ nodal
domains define a $k$-partition of $ \Omega$, with $k=\mu(u)$. 

We  recall that the  Courant nodal
 Theorem \cite{Cou}  (1923) says that, for  $k\geq 1$,  and if $\lambda_k$ denotes  the
$k$-th eigenvalue  of
the Dirichlet Laplacian in $\Omega$  and  $ E(\lambda_k)$ the eigenspace associated with $\lambda_k\,$, then,   for all real-valued $ u\in  E(\lambda_k)\setminus \{0\}\;,\;
\mu (u)\le k\,.
$
In dimension $1$ the Sturm-Liouville theory says that we have
always equality (for Dirichlet in a bounded interval) in the previous theorem (this is what we will call
later a Courant-sharp situation). 
 A  theorem due to Pleijel \cite{Pl} in
1956  says
that this cannot be true when the dimension (here we consider the
$2D$-case) is larger than one. The proof involves lower bounds for the energy of nodal partitions but what is only used is actually that the ground state energy in each  of the domain of the partition is the same.  This is this link between the proof of Pleijel's Theorem and the lower bounds for the energy of a partition that we would like to explore in this survey, motivated by two recent contributions of J. Bourgain and S. Steinerberger.
We will  focus on the large $k$ problem (and the hexagonal conjecture). This leads us also to discuss  some  conjecture by I.~Polterovich \cite{Pol}, which could be the consequence of a "square" conjecture for  some  still unknown subclass of partitions. We finally discuss some new consequence of the magnetic characterization of minimal partitions \cite{HHO3}  for the critical points of this partition.

\section{A reminder on  minimal spectral partitions}

We now introduce for $k\in \mathbb N$ ($k\geq 1$),
 the notion of $k$-partition. We
call {\bf  $ k$-partition}  of $ \Omega$ a family 
$ \mathcal D=\{D_i\}_{i=1}^k$ of mutually disjoint sets in $\Omega$.
We  call it {\bf open} if the $D_i$ are open sets of
$ \Omega$,
{\bf connected}
 if the $ D_i$ are connected. We denote by $ \mathfrak O_k(\Omega)$ the set of open connected
partitions of $\Omega$. 
We now introduce the notion of spectral minimal partition sequence.
\begin{definition}\label{regOm}~\\
For any integer $ k\ge 1$,
 and for $ \mathcal D$ in $ \mathfrak O_k(\Omega)$, we
introduce 
\begin{equation}\label{LaD}
\Lambda(\mathcal D)=\max_{i}\la(D_i).
\end{equation} 
Then we define
\begin{equation}\label{frakL} 
\mathfrak L_{k}(\Omega)=\inf_{\mathcal D\in \mathfrak O_k}\:\Lambda(\mathcal D).
\end{equation}
and  call  $ \mathcal D\in \mathfrak O_k$ a minimal $k$-partition if 
  $ \mathfrak L_{k}=\Lambda(\mathcal D)$. 
\end{definition}
More generally we can define, for $p\in [1,+\infty)$,  $ \Lambda^p(\mathcal D)$ and $\mathfrak L_{k,p}(\Omega)$ by:
\begin{equation}\label{frakLp} 
\Lambda^p(\mathcal D):=\left(\frac{\sum \lambda (D_i)^p}{k}\right)^{\frac 1p}\;,\; \mathfrak L_{k,p}(\Omega)=\inf_{\mathcal D\in \mathfrak O_k}\:\Lambda^p(\mathcal D).
\end{equation}
Note that we can minimize over non necessarily connected partitions and get the connectedness of the minimal partitions as a property (see \cite{HHOT}).

If $ k=2$, it is rather well known  that  $ \mathfrak L_2 =\lambda_2$ and
 that the associated minimal  $ 2$-partition is
 a {\bf nodal partition}, i.e. a partition whose elements
 are  the nodal domains of
   some eigenfunction corresponding to
 $ \lambda_2$. 

A partition $ \mathcal D=\{D_i\}_{i=1}^k$ of  $
\Omega$ in $ \mathfrak
O_k$ is called {\bf strong} if 
\begin{equation}\label{defstr}
\Inte(\overline{\cup_i D_i}) \setminus \pa \Om =\Om\;.
\end{equation}

Attached to a strong  partition, we  associate a closed
set in $ \overline{\Omega}$, which is called the {\bf boundary set}  of the partition~:
\begin{equation}\label{assclset} 
N(\mathcal D) = \overline{ \cup_i \left( \partial D_i \cap \Omega
  \right)}\;.
\end{equation}
$ N(\mathcal D)$ plays the role
 of the nodal set (in the case of a nodal partition).

This suggests the following definition of regularity for a partition: 
\begin{definition}\label{AMS}~\\
We call a partition  $\mathcal D$ regular if its associated
 boundary set  $ N(\mathcal D) $, has the following properties~:\\
(i)
Except for finitely many distinct $ x_i\in \Om\cap N$
 in the neighborhood of which $ N$ is the union of $\nu_i= \nu(x_i)$
smooth curves ($ \nu_i\geq 3$) with one end at $ x_i$,  
$ N$ is locally diffeomorphic to a regular 
curve.\\
(ii)
$ \pa\Om\cap N$ consists of a (possibly empty) finite set
of points $ z_i$. Moreover  
 $N$ is near $ z_i$ the union 
of $ \rho_i$ distinct smooth half-curves which hit
$ z_i$.\\
(iii) $ N$  has the {\bf equal angle
  meeting
 property}
\end{definition}
The $x_i$ are called the critical points and define the set
$X(N)$.  A particular role is played by $X^{odd}(N)$ corresponding to the critical points for which $\nu_i$ is odd.  
By {\bf equal angle meeting property}, we mean that   the half curves meet with equal angle at each critical
 point of $ N$ and also at the boundary together with the
 tangent to the boundary.

 We say that two elements of the partition $ D_i,D_j$ are {\bf  neighbors}
and write $ D_i\sim D_j$,  if $
D_{ij}:=\Inte(\overline {D_i\cup D_j})\setminus \pa \Om$ is
connected. 
We associate with  
each $ \mathcal D$ a {\bf graph}
  $ G(\mathcal D)$ by
associating with  each $ D_i$ a vertex and to each 
pair $ D_i\sim D_j$ an edge. We will say that the graph is
{\bf bipartite} if it
can be colored by two colors (two neighbors having two different
colors). We recall that the graph associated
 with  a collection of nodal domains of an eigenfunction is always
 bipartite.

\section{Pleijel's Theorem revisited}
Pleijel's Theorem as stated in the introduction is the consequence of a more precise theorem and  the aim of this section is to present a formalized proof of the historical statement permitting  to understand recent improvements and to formulate conjectures.\\

Generally, the classical proof is going through the proposition
\begin{proposition}\label{Prop1}
\begin{equation}\label{bourgain3} 
\limsup_{n\ar +\infty} \frac{\mu(\phi_n)}{n}\leq \frac{4 \pi }{ A(\Omega)  \liminf_{k\ar +\infty} \frac{\mathfrak L_k(\Omega)}{k}} \,.
\end{equation}
\end{proposition}
Here $\mu(\phi_n)$ is the cardinal of the nodal components of $\Omega \setminus N(\phi_n)$ and $A(\Omega)$ denotes the area of $\Omega$\,.  \\
Then one  establishes  a lower bound for  $ A(\Omega) \liminf_{k\ar +\infty} \frac{\mathfrak L_k(\Omega)}{k}\,$, which should  be $> 4\pi$. This property is deduced in \cite{Pl} from the Faber-Krahn inequality which says:
\begin{equation}\label{FK}
{\rm (Faber-Krahn)}\quad A(D) \lambda (D) \geq \lambda(Disk_1)\,,
\end{equation}
for any open set $D$. Here $Disk_1$ denotes the disk of area $1$.\\
Behind the statement of Proposition \ref{Prop1}, we have actually the stronger proposition:
\begin{proposition}\label{Prop2}
\begin{equation}\label{bourgain3b} 
\limsup_{n\ar +\infty} \frac{\mu(\phi_n)}{n}\leq \frac{4 \pi }{ A(\Omega)  \liminf_{k\ar +\infty} \frac{L_k(\Omega)}{k}}\,.
\end{equation}
\end{proposition}
Here $L_k(\Omega)$ denotes the smallest eigenvalue (if any)  for which there exists an eigenfunction in $E(L_k)$ with $k$ nodal domains. If no eigenvalue has this property, we simply write $L_k(\Omega) = +\infty$.\\
The proof of Proposition \ref{Prop2}  is immediate observing first  that for any subsequence $n_\ell$, we have
$$
\frac{\lambda_{n_\ell}} {n_\ell} \geq \frac{L_{\mu(\phi_{n_\ell})}} {n_\ell} =  \frac{L_{\mu(\phi_{n_\ell})}}   {\mu(\phi_{n_\ell})} \cdot \frac {\mu(\phi_{n_\ell}) } {n_\ell} \,.
$$
If we choose the subsequence $n_\ell$ such that 
$$
\lim_{\ell \ar +\infty} \frac {\mu(\phi_{n_\ell}) } {n_\ell} = \limsup_{n\ar +\infty} \frac{\mu(\phi_n)}{n}\,,
$$
we observe that, by Weyl's formula, we have:
\begin{equation} \label{Wf}
N(\lambda) \sim \frac{A(\Omega)}{4\pi}  \lambda\,,
\end{equation}
which implies
$$
\lim_{\ell \ar + \infty} \frac {\lambda_{n_\ell} }  {n_\ell} =4 \pi/A(\Omega)\,.
$$
We also have
$$
\liminf_{\ell \ar +\infty}  \frac{L_{\mu(\phi_{n_\ell})}} {\mu(\phi_{n_\ell)} } \geq \liminf_{k\ar +\infty} \frac{L_k(\Omega)}{k}\,.
$$
Hence we get the proposition.
$\square$

Proposition \ref{Prop1} is deduced from Proposition \ref{Prop2} by observing that it was established in \cite{HHOT} that
\begin{equation}
\lambda_k(\Omega)  \leq \mathfrak L_k(\Omega)  \leq L_k(\Omega)\,.
\end{equation}
The left hand side inequality is a consequence of the variational characterization of $\lambda_k$ and the right hand side is an immediate consequence of the definitions.
Moreover, and this is a much deeper theorem of \cite{HHOT}, the equalities  $\mathfrak L_k (\Omega) = L_k(\Omega)$ or $\mathfrak L_k (\Omega) = \lambda_k(\Omega)$ imply $\mathfrak L_k (\Omega) = L_k(\Omega)=\lambda_k(\Omega)$. We say that, in this case, we are in a Courant sharp situation.

If we think that only nodal partitions are involved in Pleijel's Theorem, it could be natural to introduce $\mathfrak L_k^{\sharp}(\Omega)$  where we take the infimum over
 a smaller non-empty class  of  $k$-partitions $\mathcal D =(D_1,\cdots, D_k)$. We call $\mathcal O_k^{\sharp}$ this yet undefined  class, 
which should contain all the nodal $k$-partitions, if any.  Natural  candidates for  $\mathcal O_k^{\sharp}$  will be discussed in Section  \ref{s6}.
 \begin{definition}\label{def3.3}
\begin{equation}
\mathfrak L_k^{\sharp}(\Omega) := \inf _{\mathcal D\in \mathcal O_k^{\sharp}} \max \lambda (D_i)\,.
\end{equation}
\end{definition}
Of course we have always
\begin{equation}\label{zz1}
\lambda_k(\Omega) \leq \mathfrak L_k(\Omega) \leq \mathfrak L_k^{\sharp}(\Omega) \leq L_k (\Omega)\,.
\end{equation}
Hence we have:
\begin{proposition}\label{thirdprop}
\begin{equation}\label{bourgain3c} 
\limsup_{n\ar +\infty} \frac{\mu(\phi_n)}{n}\leq \frac{4 \pi }{ A(\Omega)  \liminf_{k\ar +\infty} \frac{\mathfrak L_k^{\sharp}(\Omega)}{k}} \,,
\end{equation}
\end{proposition}
Hence this is the right hand side of \eqref{bourgain3c} which seems to be interesting to analyze.\\
It is clear from \eqref{zz1} that all these upper bounds are less than one, which corresponds to a weak asymptotic version of Courant's Theorem.\\

We now come back to the proof by Pleijel of his theorem. We apply the Faber-Krahn inequality \eqref{FK}  to any element $D_i$ of the minimal $k$-partition $\mathcal D$, and summing up,  we immediately get:
\begin{equation}\label{FKLk}
A(\Omega) \frac{\mathfrak L_k(\Omega)}{k}  \geq \lambda(Disk_1)\,.
\end{equation}
Implementing this inequality in Proposition \ref{Prop2}, we immediately get:
\begin{theorem}[Pleijel]
\begin{equation}\label{nupl}
 \limsup_{n\ar +\infty} \frac{\mu(\phi_n)}{n} \leq \nu_{Pl}\,,
\end{equation}
with $$\nu_{Pl}= \frac{4\pi}{\lambda(Disk_1)}\sim 0.691\,.$$
\end{theorem}
\begin{remarks}~
\begin{enumerate}
\item 
We note that the proof of Pleijel uses only a weak form of \eqref{FKLk}, where $\mathfrak L_k$ is replaced by $L_k$. 
\item Note that the same result is true in the Neumann case (Polterovich \cite{Pol}) under some analyticity assumption on the boundary. Note also that computations for the square were already presented in \cite{Pl}.
\item Note that  we have the better:
\begin{equation}\label{lbabc}
A(\Omega)  \frac{\mathfrak L_{k,1}(\Omega)}{k}  \geq \lambda(Disk_1)\,.
\end{equation}
But this improvement has no effect on Pleijel's Theorem. In particular,  we recall that we do not have necessarily $\lambda_k(\Omega) \leq \mathfrak L_{k,1}(\Omega)$ (take $k=2$ and use the criterion of Helffer--Hoffmann-Ostenhof \cite{HHO}).  This inequality can be replaced (see \cite{HHO}) by:
\begin{equation}
\mathfrak L_{k,1} (\Omega) \geq \frac 1k \sum_{\ell=1}^k \lambda_\ell(\Omega)\,.
\end{equation}
Again a Weyl asymptotic shows that this last inequality is strict for $k$ large. We have indeed as $k\ar +\infty$
$$
\frac 1k \sum_{\ell=1}^k \lambda_\ell(\Omega) \sim \frac{2 \pi k}{ A(\Omega)}\,,
$$ to compare with \eqref{lbabc}.

\item Pleijel's Theorem is valid in the case of the Laplace-Beltrami operator on a compact manifold (see some survey inside  \cite{BH} or \cite{Be}).
\end{enumerate}
\end{remarks}
\paragraph{Around the Hexagonal conjecture}~\\ 
It is rather easy (see \cite{BHV}) using hexagonal tilings  to prove that:
\begin{equation}\label{majhexa}
A(\Omega) \liminf_{k\ar +\infty} \frac{\mathfrak L_k(\Omega)}{k}\leq A(\Omega) \limsup_{k\ar +\infty} \frac{\mathfrak L_k(\Omega)}{k}  \leq \lambda(Hexa_1)\,,
\end{equation}
where $Hexa_1$ is the regular hexagon of area $1$.\\
Note that any tiling leads to a similar upper bound but $\lambda(Hexa_1)$ gives to our knowledge the smallest lower bound for a fundamental cell of area~$1$. 
Here are a few numerical (sometimes exact) values corresponding to the 
$Hexa_{1}$, $T_{1}$, and  $Sq_{1}$ being respectively a regular hexagon,  a square of area $1$ and an equilateral triangle:
\begin{equation}
\lambda_{1}(Hexa_{1})\sim 18.5901\;,\;
\lambda_{1}(Sq_{1})=2\pi^2\sim 19.7392 \;,\; \lambda_{1}(T_{1})\sim 22.7929\,.
\end{equation}
In addition it is not known that the regular hexagon with area 
$1$ has the lowest
groundstate eigenvalue among all hexagons of the same area.  (famous 
conjecture of Polya and Szeg\"o).

A now well known conjecture (hexagonal conjecture) (Van den Berg, Caffarelli-Lin \cite{CL}) was discussed in Helffer--Hoffmann-Ostenhof--Terracini \cite{HHOT}, Bonnaillie-Helffer-Vial \cite{BHV},  Bourdin-Bucur-Oudet  \cite{BBO} and reads as follows:
\begin{conjecture}
\begin{equation}
A(\Omega) \liminf_{k\ar +\infty} \frac{\mathfrak L_k(\Omega)}{k} = A(\Omega) \limsup_{k\ar +\infty} \frac{\mathfrak L_k(\Omega)}{k}  = \lambda(Hexa_1)
\end{equation}
\end{conjecture}
 
The minimal partitions corresponding to $\mathfrak L_{k,1}$ were computed for the torus by Bourdin-Bucur-Oudet \cite{BBO}.
This conjecture would lead to the conjecture that in Pleijel's estimate we have actually:
\begin{conjecture}{(Hexagonal conjecture for Pleijel)}
\begin{equation}
A(\Omega) \limsup_{n\ar +\infty} \frac{\mu(\phi_n)}{n} \leq \nu_{Hex}\,,
\end{equation}
with $$\nu_{Hex }= \frac{4\pi}{\lambda(Hexa_1)}\sim 0.677\,.$$
\end{conjecture}
We note indeed
that
$$
\frac{\nu_{Hex }}{\nu_{Pl}} = \frac{ \lambda(Disk_1)}{\lambda(Hexa_1)}\sim 0.977\,.
$$

We now come back to the enigmatic Proposition  \ref{thirdprop} and consider the asymptotic behavior of $\frac{\mathfrak L_k^{\sharp}(\Omega)}{k}$ as $k\ar + \infty$. A first 
indication that our choice of
$\Om^\#$ is promising would be to show that the following property holds.
 \begin{property}\label{property}
\begin{equation}
A(\Omega) \liminf _{k\ar +\infty} \frac{\mathfrak L_k^{\sharp}(\Omega)}{k} \geq  \lambda (Sq_1)\,.
\end{equation}
where $Sq_1$ denotes the unit square.\\
\end{property}

The proof of this  property should mimic  what was done for  $\mathfrak L_k(\Omega)$ (see for example \cite{BHV} or \cite{CL}), but replacing  hexagonal tilings by square tilings. 
 \begin{figure}[h!]\label{hexasquare}
\begin{center}
\includegraphics[width=7cm]{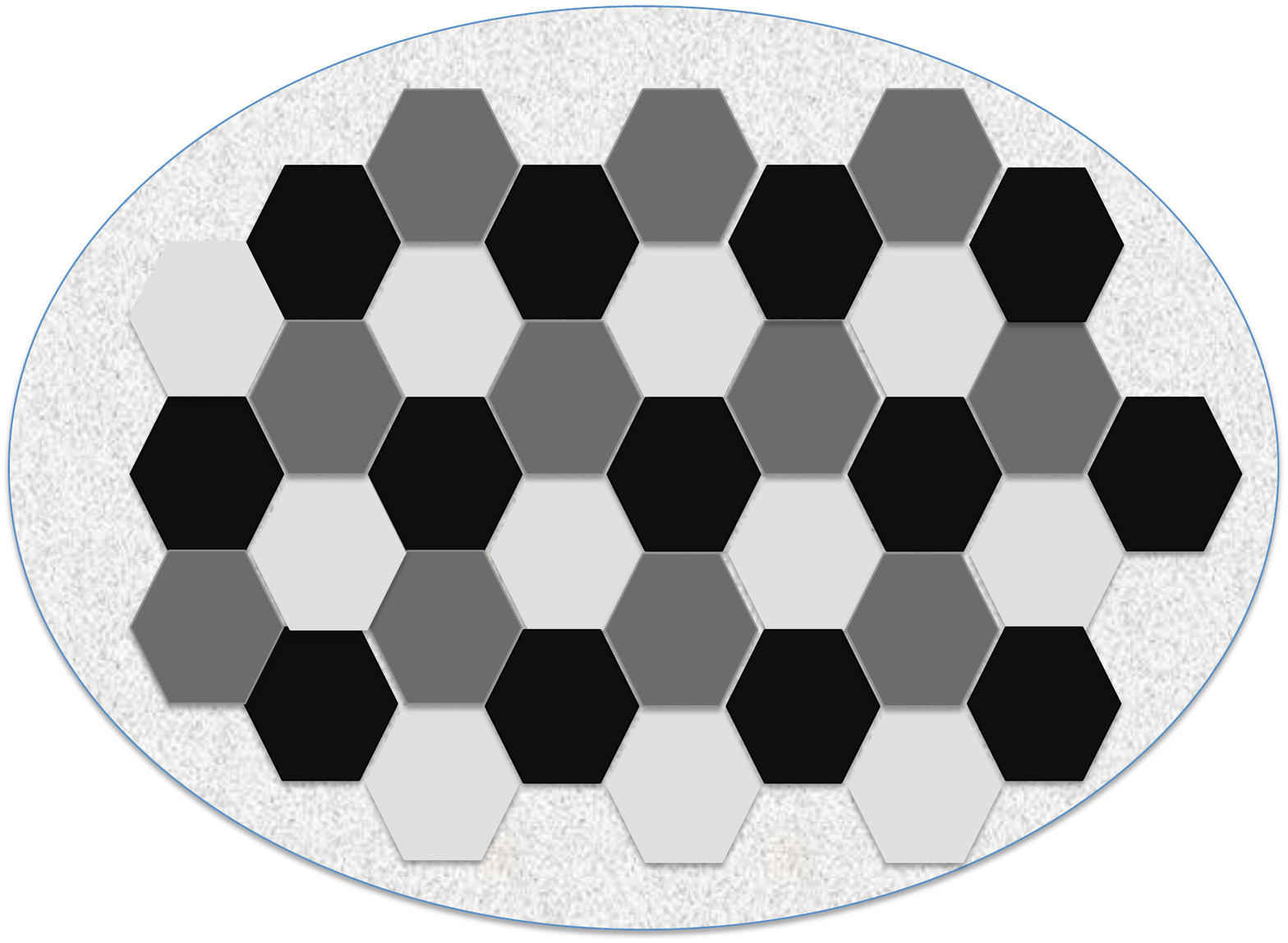}~\includegraphics[width=6cm]{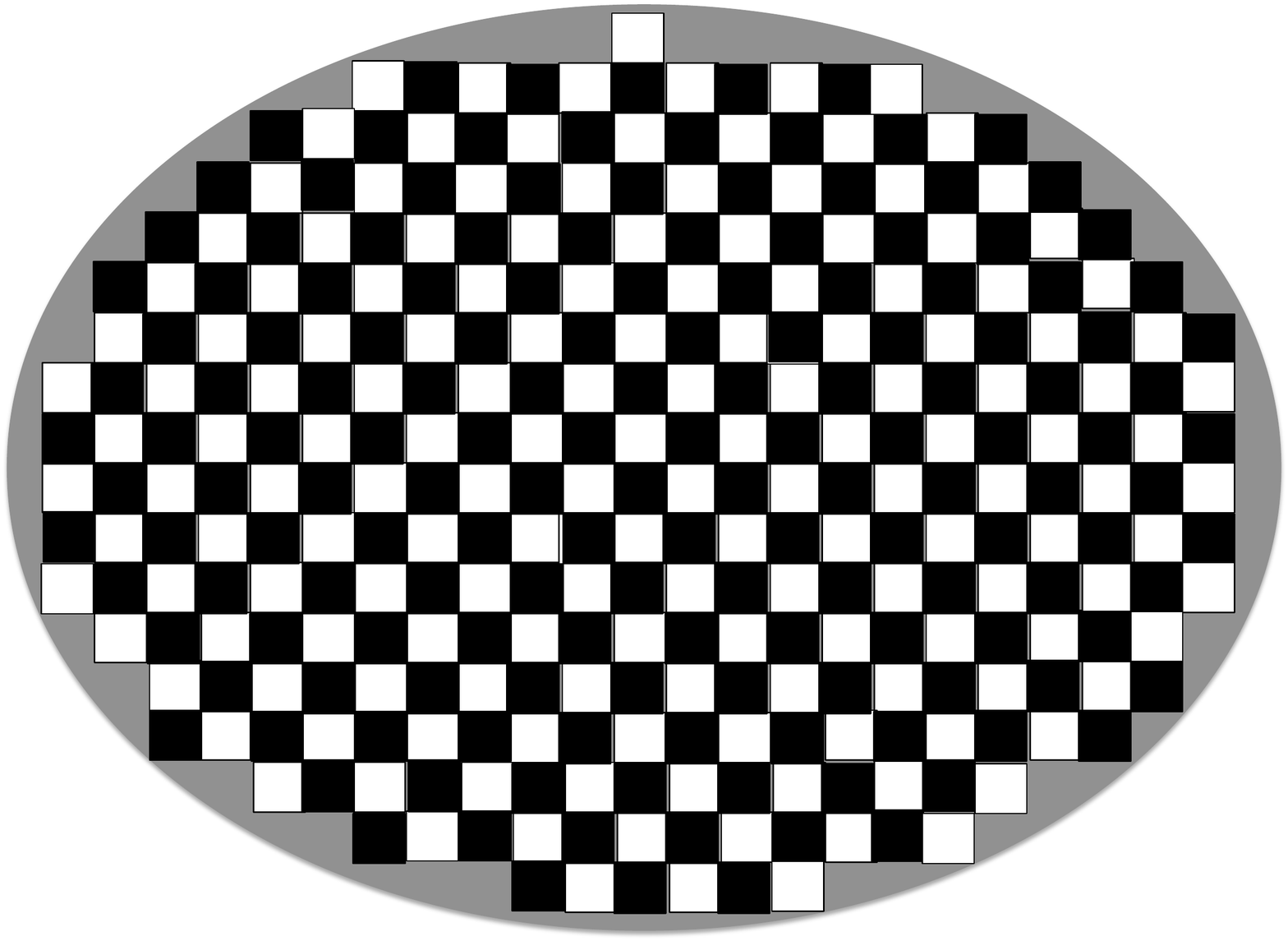}
\caption{Non exhaustive hexagonal or square tilings inside an open set $\Omega$}
 \end{center}
 \end{figure}

The philosophy behind the choice of $\Omega_k^\sharp$ should be the following: the hexagonal conjecture for $k$-partitions should be replaced by the square conjecture when bipartite $k$-partitions are involved because square tilings can be colored by two colors with the rule that two neighbors have two different  colors.

Note that the  existence of classes $\mathcal O_k^\sharp$ ($k\in \mathbb N^*$) such that Property \ref{property} is satisfied would give a proof of  the conjecture:
\begin{conjecture}{(Polterovich)}
\begin{equation}\label{bourgain3e} 
\limsup_{n\ar +\infty} \frac{\mu(\phi_n)}{n}\leq \frac{4 \pi }{\lambda(Sq_1)}=\frac{2}{\pi}\,. 
\end{equation}
\end{conjecture}
This conjecture is due to Iosif Polterovich \cite{Pol} on the basis of computations for the rectangle of Blum-Gutzman-Smilansky \cite{BGS}. Due to the computations on the square \cite{Pl} (see however the discussion in Section \ref{Sectionrectangle}), together with computations for  the rectangle \cite{SS}, this should be the optimal result.

\section{Improving the use of the Faber-Krahn Inequality}
\subsection{Preliminaries}
This section is devoted to reporting on the two previously mentioned results by J. Bourgain and S. Steinerberger. Although not explicitly written in this way, the goal of Bourgain  \cite{Bo} and Steinerberger \cite{St} was to improve Pleijel's proof by improving the lower bound of $ \liminf_{k\ar +\infty} \frac{\mathfrak L_k(\Omega)}{k}$.  Bourgain gives a  rough  estimate of his improvement  with a  size of $10^{-9}$.\\
 
 In any case, it is clear from their proof  that  this will lead to a statement where in  \eqref{nupl} $\nu_{Pl}$ is replaced by 
 $$
 \nu_{Hex}\leq \nu_{Bo} < \nu_{Pl}\,,
 $$
 and
 $$
 \nu_{Hex} \leq \nu_{St} < \nu_{Pl}\,,
 $$
where $\nu_{Bo}$ and  $\nu_{St}$ are the constants obtained respectively by  Bourgain and Steinerberger.\\

\subsection{Bourgain's improvement \cite{Bo}}
One ingredient is a refinement of the Faber-Krahn inequality du to Hansen-Nadirashvili \cite{HaNa}:
\begin{lemma}[Hansen-Nadirashvili]~\\
 For a nonempty simply connected bounded domain $ \Omega \subset  \mathbb R^2$, we have
\begin{equation}\label{NH}
A(\Omega)\, \lambda(\Omega) \geq  \left(1 + \frac{1}{250} (1- \frac{r_i(\Omega)}{r_0(\Omega)})^2 \right)\lambda ({\rm Disk_1} )\,,
 \end{equation}
with $  r_0(\Omega)$ the radius of the disk of same area as $ \Omega$  and $ r_i(\Omega)$  the inradius of $ \Omega$.
 \end{lemma}
 
 Actually, J. Bourgain  needs (and gives) an additional argument for treating non simply connected domains. 
In the right hand side of \eqref{NH}  not only the inradius occurs but also the smallest area of the 
components of $\mathbb R^2\setminus \Om$.\\
 
 The other very tricky idea  is to use quantitatively that all the open sets of the partition cannot be very close to disks (packing density) (see Blind \cite{Blind}).

The inequality obtained by Bourgain is the following (see (37) in \cite{Bo}) as $ k\ar +\infty$, is that for any $ \delta\in (0,\delta_0)$
\begin{equation}\label{bourqain}  
\frac{\mathfrak L_k (\Omega)}{k} \geq (1 + o(1)) \lambda({\rm Disk_1} )A(\Omega)^{-1} \times b(\delta) 
\end{equation}
where
$$
 b(\delta):= ( 1 + 250 \delta^{-3}) (\frac{\pi}{\sqrt{12}}(1-\delta)^{-2} + 250 \delta^{-3})^{-1}\,.
$$
and $ \delta_0 \in (0,1)$ is computed with the help of the packing condition. This condition reads
$$ 
\frac{\delta_0^3}{250} = (\frac{1-\delta_0}{p})^2 -1\,,
$$
where $ p$ is a packing constant determined by Blind \cite{Blind} ($  p \sim 0.743\,$).

But for $ \delta>0$ small enough, we get $ b(\delta) >1$ (as a consequence of $ \frac{\pi}{\sqrt{12}} <1$), hence Bourgain  has improved what was obtained via Faber-Krahn (see \eqref{FKLk}).\\
As also observed by Steinerberger, one gets
$$
 \frac{\lambda({\rm Hexa_1} ) }{\lambda({\rm Disk_1} )} \geq \sup_{\delta\in(0,\delta_0)} b(\delta)>1\,,
$$
which gives a limit for any improvement of the estimate.\\
In any case, we have
\begin{equation}\label{bourgain1}
 \liminf_{k \ar + \infty}  \frac{\mathfrak L_k (\Omega)}{k} \geq  \lambda({\rm Disk_1} )A(\Omega)^{-1} \times \sup_{\delta \in (0,\delta_0)} b(\delta) \,.
\end{equation}

\subsection{The uncertainty principle by S. Steinerberger}
To explain this principle, we associate with  a partition $\mathcal D =( \Omega_i)$ of $ \Omega$
$$ 
D(\Omega_i) = 1 - \frac{\min_j A(\Omega_j)}{A(\Omega_i)}\,.
$$
We also need to define for an open set  $D$ with finite area, the Fraenkel asymmetry of $D$:
$$
\mathcal A(D) = \inf_B \frac{A(D\triangle B)}{A(D)}\,,
$$
where the infimum is over the balls of same area and where  $$D \triangle B = (D \setminus B) \cup (B\setminus  D)\,.$$
Steinerberger's uncertainty principle reads:
\begin{theorem}~\\
There exists a universal constant $ c>0$, and a  $ k_0(\Omega)$ such that for each $k$-partition of $\Omega$: $\mathcal D =(\Omega_i)_{i=1,\dots,k}\,,$ with $k \geq k_0(\Omega)$, 
\begin{equation}\label{UP}
\sum_i (D(\Omega_i) +\mathcal A(\Omega_i)) \,  \frac{A(\Omega_i)}{A(\Omega)}  \geq c\,.
\end{equation}
\end{theorem}

\subsection{Application to equipartitions of energy $ \lambda$}
Let us show how we recover a lower bound for $ \liminf_{k\ar +\infty} \left( \mathfrak L_k(\Omega)/k \right)$ improving \eqref{FKLk} asymptotically. We consider a $ k$-equipartition of energy $ \lambda$.  We recall from \cite{BH} that an equipartition is a strong partition for which the ground state energy in each open set $D_i$ is the same.  In particular, nodal partitions and minimal partitions for $\mathfrak L_k$  are typical examples of equipartitions. If we assume that $k \geq k_0(\Omega)$, the uncertainty principle says that its is enough to consider two cases.\\
We first assume that 
$$
\sum_i D(\Omega_i) \frac{A(\Omega_i)}{A(\Omega)}  \geq \frac c2 \,.
$$
We can rewrite this inequality in the form:
$$
 k \inf_j A(\Omega_j) \leq (1 -\frac c2) A(\Omega)\,.
$$
After implementation of Faber-Krahn, we obtain
\begin{equation}\label{case1} 
\frac{k}{\lambda} \lambda(Disk_1) \leq (1- \frac c 2) A(\Omega)\,.
\end{equation}
We now assume that 
$$
\sum_i \mathcal A(\Omega_i)  \frac{A(\Omega_i)}{A(\Omega)}  \geq \frac c2 \,.
$$
This assumption implies, using that $\mathcal A(\Omega_i)\leq 2\,$, 
\begin{equation}\label{info}  
A \left( \cup_{\{\mathcal A(\Omega_i) \geq \frac c6\}} \Omega_i \right) \geq \frac c 6 A(\Omega)\,.
\end{equation}
The role of $ \mathcal A$ can be understood from the following inequality due to Brasco, De Philippis, and Velichkov \cite{BDPV}:\\
There exists $ C>0$ such that, for any open set  $ \omega$ with finite area, 
\begin{equation}\label{ineqbpv}
 A(\omega) \lambda(\omega) - \lambda (Disk_1)\geq C \mathcal A(\omega)^2  \lambda(Disk_1)\,.
\end{equation}
If we apply this inequality with $ \omega =\Omega_i\,$, 
it reads
$$ A(\Omega_i) \lambda - \lambda (Disk_1)\geq C \mathcal A(\Omega_i)^2 \lambda(Disk_1)\,.
$$
 Hence we get for any $ i$ such that $ \mathcal A(\Omega_i)\geq \frac c6\,$,  the inequality 
\begin{equation}\label{case2} 
\lambda(Disk_1) (1 + \frac{Cc^2}{3 6})  \leq  A(\Omega_i) \lambda \,,
\end{equation}
which is an improvement of Faber-Krahn for these $ \Omega_i$'s.\\
Summing over $ i$  and using  \eqref{info} leads to
\begin{equation*}\label{case20} 
\frac{k}{\lambda} \lambda(Disk_1) \leq  (1 + \frac{Cc^2}{3 6})^{-1}  \,   A(\Omega) \left (1 + (1-\frac c6) \frac{Cc^2}{36}\right)\,,
\end{equation*}
and finally to
\begin{equation}\label{case21} 
\frac{k}{\lambda} \lambda(Disk_1) \leq   \, 
 \left (1 - \frac{Cc^3}{216+ 6 C c^2} \right)
A(\Omega)\,. \end{equation}

Putting \eqref{case1} and \eqref{case2} together, we obtain that for $ k\geq k_0(\Omega)$ (as assumed from the beginning) the $ k$-partition satisfies
\begin{equation}\label{case1+2} 
\frac{k}{\lambda} \lambda(Disk_1) \leq  \max \left(  (1-\frac c2),  (1 - \frac{Cc^3}{216+ 6 C c^2} ) \right)   A(\Omega)\,.
\end{equation}
If we apply this to minimal partitions ($ \lambda =\mathfrak L_k(\Omega)$), this reads
\begin{equation}\label{final} 
\lambda(Disk_1) \leq  \max \left(  (1-\frac c2), (1 - \frac{Cc^3}{216+ 6 C c^2} ) \right)   A(\Omega) \liminf_{k\ar +\infty} \frac{\mathfrak L_k(\Omega)}{k}\,.
\end{equation}
Hence S. Steinerberger  recovers Bourgain's improvement \eqref{bourgain1} with a not explicit constant\footnote{At least $C$ in \eqref{ineqbpv}  is not explicit.}.

\begin{remark}
Steinerberger obtains also a similar lower bound to \eqref{final}  for $\liminf_{k\ar +\infty} \frac{\mathfrak L_{k,1}(\Omega)}{k}$ using a convexity argument.
\end{remark}
\section{Considerations around rectangles} \label{Sectionrectangle}
The detailed analysis of the spectrum of the Dirichlet Laplacian in a rectangle is basic in Pleijel's paper \cite{Pl}. As mentioned in \cite{CH}, other results have been previously obtained in the PHD of A. Stern \cite{Stern}, defended in 1924. Other aspects relative to spectral minimal partitions appear in \cite{HHOT}. Take $\mathcal R(a,b)=(0,a\pi)\times (0,b\pi)$. The eigenvalues are given by 
$$\hat \la_{m,n}=(\frac{m^2}{a^2} +\frac{n^2}{b^2})\,,$$
with a corresponding basis of eigenfunctions given by
$$
\phi_{m,n} (x,y)= \sin \frac{mx}{a} \,\sin \frac{ny}{b}\,.
$$
If it is easy to determine the Courant sharp cases when $b^2/a^2$ is irrational (see for example  \cite{HHOT}). The rational case is more difficult. Pleijel claims in \cite{Pl} that in the case of the square it is Courant sharp if and only if  $k=1,2,4$.  The exclusion of $k=5, 7,9$ is however not justified (the author refers indeed to Courant-Hilbert \cite{CH} where only pictures are presented, actually taken from the old book (1891) by Pockels \cite{Poc}). This can actually be controlled by an explicit computation of the nodal sets of each combination  ($\theta \in [0,\pi))$:  $$ \Phi_{m,n,\theta}:=\cos \theta \phi_{m,n} + \sin \theta \phi_{n,m}$$
for $(m,n) = (1,3)\,,\, (1,4)\,$, and \, $(2,3)$.

In this context the following guess could  be natural:\\
{\it Suppose that $\hat \la_{m,n}$ has multiplicity $\mathfrak m(m,n)\,$. 
 Let $\mu_{\rm max}(u)$ be the maximum of the number of nodal domains of the eigenfunctions in the eigenspace associated with
$\hat \la_{m,n}$. 
$$
\mu_{\rm max}=\sup_j  (m_jn_j) \,,
$$
where the sup is computed over the pairs $(m_j,n_j)$ such that
$$\hat  \lambda_{m_j,n_j}=\hat \la_{m,n}\,.$$
}
The problem  is not easy because one has to consider, in the case of degenerate eigenvalues, linear combinations of the canonical eigenfunctions 
associated with the  $\hat \la_{m,n}\,$. Actually, as stated above, the guess is wrong. As observed by  Pleijel \cite{Pl},  the eigenfunction  $\Phi_{1,3, \frac{3\pi}{4}}$ corresponding to the fifth eigenvalue has four nodal domains delimited by the two diagonals, and $\mu_{max}=3$.  More generally one can consider $u_k:= \Phi_{1,3,  \frac{3\pi}{4}}( 2^k x, 2^k y)$ to get an eigenfunction associated with the eigenvalue $\lambda_{n(k)}= \hat \lambda_{2^k,3\cdot 2^k}=10 \cdot  4^k$ with $4^{k+1}$ nodal domains. 
 The corresponding quotient $\frac{\mu(u_k)}{n(k)}$  is asymptotic to $\frac{8}{5 \pi}$. This does not  contradict the Polterovich conjecture.
Note also that  for
each number $K$, there is
an eigenfunction $u$ with $\mu(u)\ge K$. Finally let us mention that counterexamples to a similar guess  in the Neumann case can be found in \cite{Poc}.\\
 
{\bf Pleijel's constant.}\\  We consider 
for each $\Om$ and each orthonormal basis $\mathcal B_\Omega := (u_n)$  of the Dirichlet Laplacian:
\begin{equation}\label{PlOmB}
Pl(\Om,\mathcal B_\Omega)=\limsup_{n\ar +\infty} \frac{\mu^{\Om,\mathcal B}_n}{n}\,,  
\end{equation}
where $\mu_n^{\Om,\mathcal B}$ denotes the number of nodal domains of  $u_n \,$.\\
The reference to $\mathcal B$ is only needed in the case when the Laplacian has a infinite sequence of multiple eigenvalues.
We then define:
\begin{equation}\label{PlOm}
Pl(\Omega)=\sup_{\mathcal B_\Omega}Pl(\Om,\mathcal B_\Omega) \,.
\end{equation}
Now the question arises how and whether $Pl(\Om)$ depends on $\Om$. Note that $Pl(\Om) =Pl(T\Om)$ where
 $T$ denotes
 scaling or rotation, reflection, translation.\\
It is not even clear that $\mu_n^{\Om,\mathcal B}$ tends for every $\Om$ to infinity, see \cite{THO:2012}. 
The Pleijel constant could be defined as 
\begin{equation}\label{Plgeneral}
 Pl=\sup_{\Om}Pl(\Om)\,,
\end{equation}
and it is not at all clear that a maximizing pair $(\Om,\mathcal B_\Omega)$ exists. (The square or more generally rectangles might be good candidates
 as mentioned above.) It would be interesting to find those domains, for which it  is possible to 
calculate $Pl(\Om)\,$. 

We finally recall (cf  \cite{BGS} and  \cite{Pol}) that
\begin{proposition}\label{Rec}
Let us assume that $b^2/a^2$ is irrational.
\begin{equation}\label{Prec}
Pl(\mathcal R(a,b))=\frac{2}{\pi}\,.
\end{equation} 
\end{proposition}
{\bf Proof}\\
It suffices to consider $\mathcal R(\pi,b\pi)$ for irrational $b^2$. Since $b^2$ is irrational the eigenvalues are simple and 
\begin{equation}\label{mumn}
\mu(\phi_{m,n})=mn\,.
\end{equation}
Weyl asymptotics tells us that with $\lambda=\hat \lambda_{m,n}\,$:
\begin{equation}\label{Weyl}
k(m,n):= \#\{(\tilde m,\tilde n):\hat \la_{\tilde m,\tilde n}(b)<\la\}=\frac{b\pi}{4}(m^2+n^2/b^2)+o(\la)\,.
\end{equation}
We have
$$\lambda_{k(m,n) +1} = \hat \lambda_{m,n}\,.
$$
We observe that $\mu(\phi_{n,m})/ \hat k(n,m)$ is asymptotically given by
\begin{equation}\label{Pmnb}
P(m,n;b):=\frac{4mn}{\pi(m^2b+n^2/b)}\le \frac{2}{\pi}\,.
\end{equation}
 Next we take a sequence $(m_k,n_k)$ such that  $b=\lim_{k\rightarrow \infty}\frac{n_k}{m_k}$ with  $m_k\ar +\infty$. \\
 We observe that
\begin{equation}\label{Paj}
\lim_{k\ar +\infty}  P(m_k,n_k,b)=\frac{2}{\pi}\,.
\end{equation}
This proves the proposition using the sequence of eigenfunctions $\phi_{m_k,n_k}$.\\
\begin{remark}
We have consequently 
\begin{equation}
Pl \geq \frac{2}{\pi}\,,
\end{equation}
the conjecture being that one has actually the equality.\\
The case when $b^2/a^2\in \qz$ depends on the discussion at the beginning of the section. We only know that
\begin{equation}
Pl (R(a,b)) \geq \frac{2}{\pi}\,,
\end{equation}

\end{remark}

\section{Looking for a class $\mathcal O^\#$.}\label{s6}
We now start the discussion on tentative choices of the classes $\mathcal O_k^{\sharp}$ (see Definition \ref{def3.3}).
\subsection{Bipartite partitions}\label{ss61}
If we think that only nodal partitions are involved in Pleijel's theorem, it could be natural to consider as class $\mathcal O_k^{\sharp}$ 
the class $\mathcal O_k^{bp}$  of the bipartite strong regular connected $k$-partitions $\mathcal D =(D_1,\cdots, D_k)$.  Note that there is some arbitrariness in the definition  but "strong" is necessary to define a bipartite partition. 
\begin{definition}
\begin{equation}\label{kbp}
\mathfrak L_k^{bp}(\Omega) := \inf _{\mathcal D\in \mathcal O_k^{bp,str}} \max \lambda (D_i)\,.
\end{equation}
\end{definition}
Although this definition is natural, all what has been established relatively to $\mathfrak L_k(\Omega)$ is unclear or at least unproved in the case of this $\mathfrak L_k^{bp}(\Omega)$.

By definition, we know that $\mathfrak L_k^{bp}(\Omega)\leq L_k(\Omega)$.  If the inequality is strict  the infimum  cannot by 
definition correspond to a nodal partition.
 If we want this notion to be helpful for improving Pleijel's constant,  it  is natural to first ask if  
$\mathfrak L_k^{bp}(\Om)>\mathfrak L_k(\Omega)$, at least for $k$ large. However we will show
\begin{proposition}\label{prop6.2}
Suppose that $\Omega$ is simply connected. Then
\begin{equation}  
\mathfrak L_k^{bp} (\Omega)= \mathfrak L_k(\Om)\,.
\end{equation}
\end{proposition}
Hence this class which could a priori appear to be  a natural candidate for $\mathcal O^\sharp$ does not lead to any improvement of the hexagonal conjecture for Pleijel's theorem.

\subsection{Proof of Proposition \ref{prop6.2}} 
{\bf Particular case.}\\
 Suppose that $\Om\subset \rz^2$, $k\geq 2$ and  consider 
a minimal $k$-partition $\mathcal D=\{D_i,\dots, D_k\}$ which is not bipartite.
 We first prove the proposition in a particular case.
\begin{lemma}
We assume that $\Om$ is simply connected and  that 
\begin{equation}\label{dD1}
 \#\{\cup \pa D_i\}=1\,.
\end{equation}
Then there is a sequence of bipartite $k$-partitions $\widehat{\mathcal D}_k(\ep)=\{\widehat D_1(\ep),\dots,\widehat D_k(\ep) \}$ of $\Omega$
with the property that
\begin{equation}  
\Lambda(\widehat{\mathcal D}_k(\ep)) \rightarrow \mathfrak L_k(\Om),\: j=1,\dots, k.
\end{equation}
\end{lemma}
For  $\ep>0$, we define for any element of the partition
\begin{equation}\label{Diep}
 D_i(\ep) =\{x\in D_i\:|\: {\rm dist}(x, \pa D_i)>\ep\}.
\end{equation}
For $\epsilon >0$ small enough all the $D_i(\epsilon)$ are non empty and connected.\\
We also define a tubular  $\epsilon$-neighborhood of $\pa \Omega \cup N(\mathcal D))$ in $\Omega$ :
$$
\mathcal S^\epsilon = \{ x\in \Omega, d(x, \pa \Omega \cup N(\mathcal D)) < \epsilon \}
$$
 $\Sg^\ep$ is connected due to Assumption \eqref{dD1}.\\
Now as $\ep$ tends to zero, $A(\Sg^\ep) \rightarrow 0$. We consider the $k$-partition $\widehat {\mathcal D} (\epsilon)$
 defined by 
 $$
 \widehat D_1(\epsilon) = D_1 \cup \mathcal S^\epsilon\,,\, \widehat D_i (\epsilon) = D_i(\epsilon)\,,\, \forall i > 1\,.
 $$
 This gives a connected open $k$-partition of $\Omega$ with the following property:
 $$
 \lambda (  \widehat D_1(\epsilon) ) < \lambda (D_1)\,,\, \lim_{\epsilon \ar 0}  \lambda (\widehat D_i (\epsilon) )= \lambda(D_i)\,,\, \forall i > 1\,.
$$
$\square$\\
Figure 2 describes the construction in the case of the disk, assuming (see \cite{HHOT}, \cite{HHO2}) that the minimal $3$-partition is the Mercedes-star.\\
\begin{figure}[h!]
\begin{center}
\includegraphics[width=6cm]{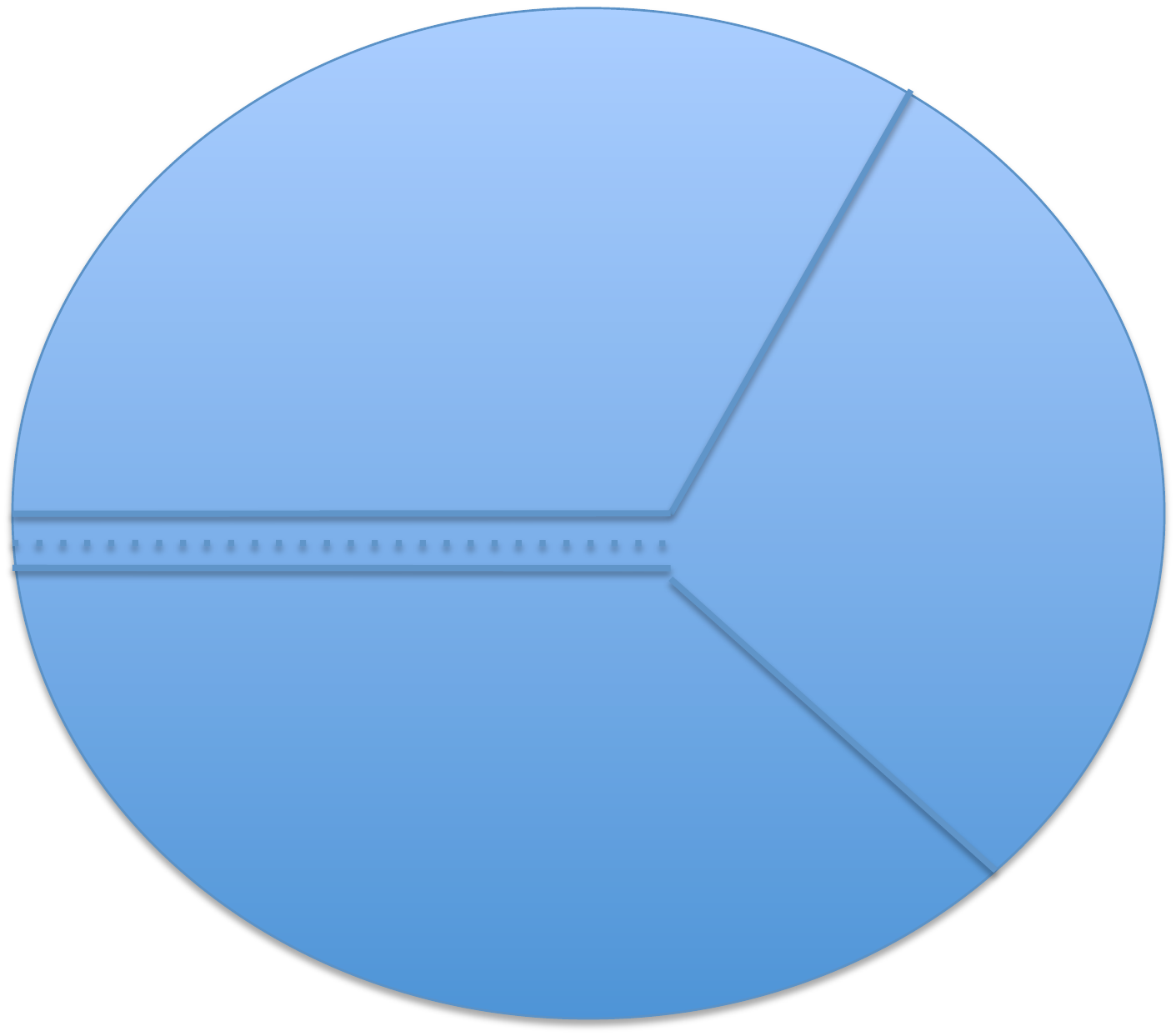}
 \end{center}
 \caption{Scheme of the construction for the Mercedes Star}
 \end{figure}

{\bf General case.}\\
  We now give the proof in the general case.  Considering the previous discussion, we can distinguish two cases for our minimal $k$-partition $\mathcal D$. 
  \begin{enumerate}
  \item $N(\mathcal D)$ does not meet $\pa \Omega$
  \item $N(\mathcal D)\cap \pa \Omega \neq \emptyset$.
  \end{enumerate}
  In the first case, after relabeling  we can call $D_1$ the unique element of the partition whose boundary touches $\pa \Omega$. We follow the previous discussion  and define $S^{(1)}(\epsilon)$ the connected component of the set $S^\epsilon$ containing $\pa D_1\cap \Omega$. The first element of the approximating $k$-partition is then $$\widehat D_1(\epsilon) :=D_1\cup S^{(1)}(\epsilon)\,.$$

  In the second case, after relabeling, we can take as $D_1$ one element of the partition such that $\pa D_1\cap \pa \Omega\neq \emptyset$ and as before introduce $\widehat D_1(\epsilon)$ as before but with  $S^{(1)}(\epsilon)$ the connected component of the set $S^\epsilon$ containing $\pa \Omega$. 
  
  We now consider  the connected components $\Omega \setminus \overline{\widehat D_1(\epsilon)}$. Many of them are simply $D_j(\epsilon)$ for $j>1$. We keep these open sets as elements of our new partition. Other components contain more than one $D_\ell$. If we denote by $\Omega^{(\ell)}(\epsilon)$ such a component, we observe that we are necessarily in case (i) of the previous discussion with $\Omega^{(\ell)}(\epsilon)$ replacing $\Omega$. Only one $D_k(\epsilon)$ inside $\Omega^{(\ell)}$ can have its boundary touching 
 $ \overline{\widehat D_1(\epsilon)}$. We can iterate inside $\Omega^{(\ell)}(\epsilon)$ what we have done in $\Omega$ and the procedure will stop after a finite number of iterations.
 \begin{remark}
The case when $\Omega$ is not simply connected can be handled similarly.
 \end{remark}
 \subsection{Almost nodal partitions}\label{s7}
Here is a new try for a definition of $\mathcal O^\sharp$  in order to have a flexible notion of partitions which are close with nodal partitions. 
 We assume that $\Omega$ is regular and simply connected.\\
 We will say that a $k$-partition $\mathcal D$  of $\Omega$  of energy $\Lambda(\mathcal D)$ is almost nodal, if there is a connected open set $\Omega'\subset \Omega$ and a $(k-1)$-subpartition $\mathcal D'$ of $\mathcal D$ such that $\mathcal D'$ is a nodal partition of $\Omega'$ of energy $\Lambda(\mathcal D)$. Of course a nodal partition is almost nodal. The first useful observation is that for any $k$ 
 there exists always an almost nodal $k$-partition.
 The proof is obtained using a sufficiently thin "square" $(k-1)$-partition in $\Omega$ and completing by the complement in $\Omega$ of the closure of the union of the preceding squares. See the right subfigure of Figure 1.
 Denoting by $\mathcal O^{a.n}_k$ the  set of the almost nodal partitions, we introduce
 \begin{equation}
 \mathfrak L_k^{a.n}(\Omega) = \inf_{\mathcal D \in \mathcal O^{a.n,j}_k} \Lambda(\mathcal D)\,.
 \end{equation}
 Of course, we have
 \begin{equation}
 \mathfrak L_k(\Omega)\leq \mathfrak L_k^{a.n}(\Omega) \leq L_k(\Omega)\,.
 \end{equation}
The next point is to observe, by the same proof giving \eqref{majhexa} but playing with the square tiling (see Figure 1), that
\begin{equation}\label{limsupsq}
A(\Omega) \limsup_{k\ar + \infty} \frac{\mathfrak L_k^{a.n}(\Omega)}{k} \leq \lambda(Sq_1)\,.
\end{equation}
Again the question arises  about the asymptotic behavior of $\liminf_{k\ar + \infty} \frac{\mathfrak L_k^{a.n}(\Omega)}{k}$. Unfortunately there are good reasons to think that we can improve \eqref{limsupsq} by proving 
 \begin{equation}
A(\Omega)  \limsup_{k\ar + \infty} \frac{\mathfrak L_k^{a.n}(\Omega)}{k} \leq \lambda(Hexa_1)\,.
 \end{equation}
 We just give an heuristical hint.  For $k$ large, we try to  "almost"  fill $\Omega$ with a maximal number of $(k-1)$ adjacent isometric regular hexagons $D_i$ ($i=1,\dots,k-1)$. 
  For $k$ large,
 they should have an area  of order $A(\Omega)/k$ and an energy of order $k \lambda(Hexa_1)/A(\Omega)$. and we complete the partition by taking as $D_k$ the complement in $\Omega$ of $\overline{\cup_{i=1}^{k-1} D_i}$. We can in addition have the property that $\lambda(D_k)\leq \lambda(D_1)$ (one way is to start with $k$ adjacent regular hexagons and to delete one). Then we construct our set $\Om'_{k-1}$ by  subtracting cracks (edges of some of the hexagons) from $\Inte (\overline{\cup_{i=1}^{k-1} D_i})$ in such a way that $(D_1,\dots, D_{k-1})$ becomes a nodal $(k-1)$-partition of $\Om'_{k-1}$ (see Figure 3). Such a construction is detailed in \cite{BHV}, when exploring the consequences of the hexagonal conjecture. This conjecture would actually impose that this is the nodal partition of a Courant sharp eigenfunction but we do not need it at this stage. Then the partition is almost nodal and asymptotically of energy $k \lambda(Hexa_1)/A(\Omega)$. \\
 \begin{figure}[h!]
\begin{center}
\includegraphics[width=6cm]{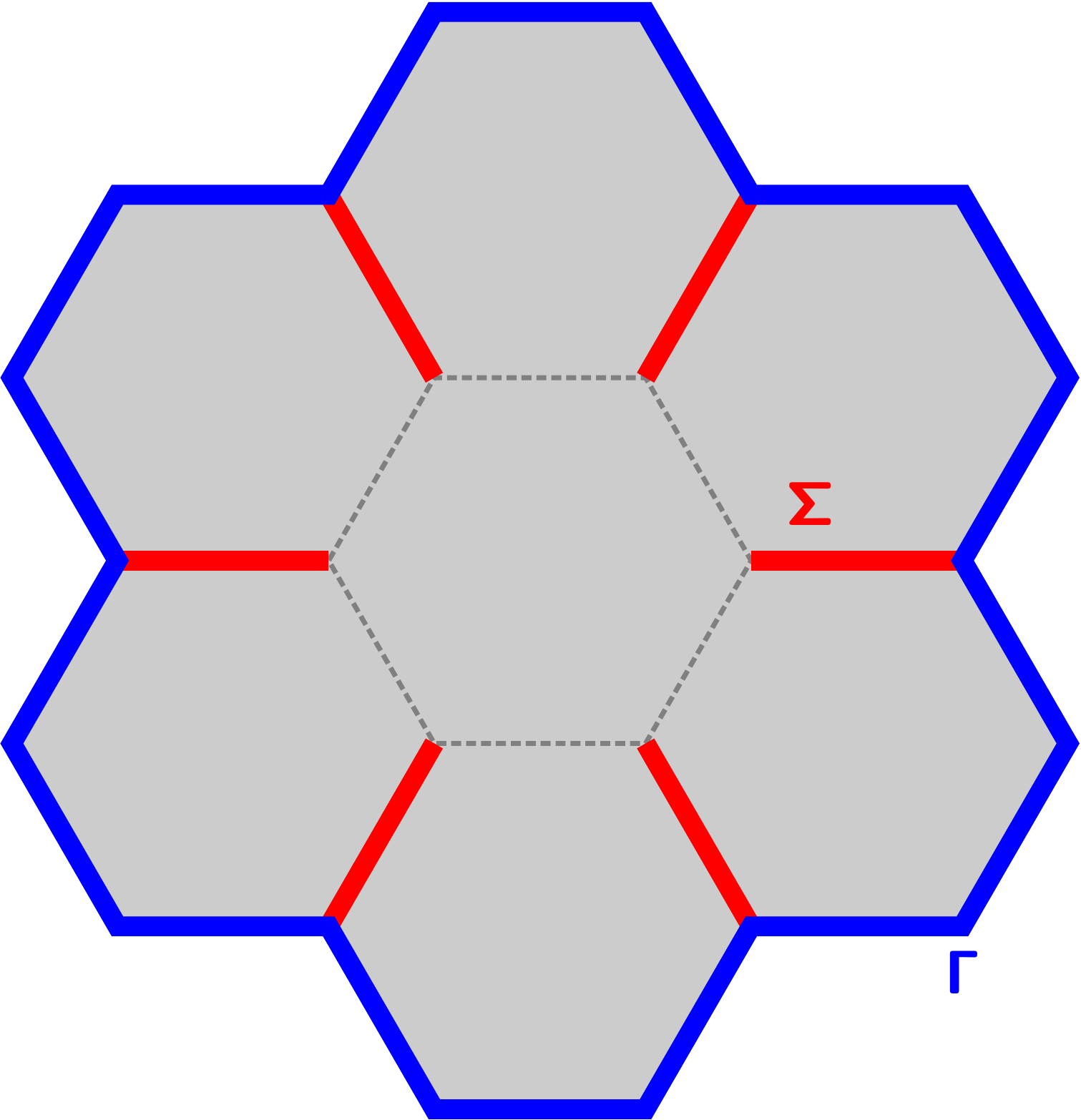}
 \end{center}
 \caption{Scheme of the construction of the cracks for $(k-1)=7$.}
 \end{figure}
Actually, starting directly from a minimal $k$-partition $\mathcal D= (D_1,\dots,D_k)$ and proceeding as before with the $(k-1)$ first elements, one can add curved segments
 belonging to the boundary of the partition such that we get a nodal partition. Here we use a property observed in \cite{HHOT} (Proof of Proposition 8.3) (see also Corollary 2.11 in \cite{BHV}). This will directly lead to the stronger
equality  $\mathfrak L_k^{a.n}(\Omega)= \mathfrak L_k(\Omega)$.\\
  Of course, one could think that by imposing more regularity on the partition and on $\Omega'$, one can eliminate this kind of examples. But as in the previous subsection, an approximation of the cracks by fine tubes could probably be used for getting the same inequality. This  we have not checked and will be  much more technical than for the proof of Proposition \ref{prop6.2}.\\
   Hence the class of almost nodal partitions is  probably too large for getting a higher infimum.
 
\subsection{Conclusion} In conclusion, we were looking for smaller classes of partitions containing nodal partitions with the hope to give some justification for the Polterovich conjecture. We have shown that two natural choices do not give a confirmation of this conjecture as initially expected.

 \section{Pleijel's Theorem for Aharonov-Bohm operators and application to minimal partitions}\label{s8}
 \subsection{The Aharonov-Bohm approach}
Let us recall some definitions and results about the Aharonov-Bohm 
Hamiltonian in an open set $\Omega$ (for short $ {\Ab\Bb}X$-Hamiltonian) with a singularity at $ X\in \Omega$ as introduced in \cite{HHOO}. We denote by $ X=(x_{0},y_{0})$ the coordinates of the pole and 
consider the magnetic potential with  flux at $ X$: 
$ \Phi = \pi  $, defined in $\dot{\Omega_X}= \Omega \setminus \{X\}$:
\begin{equation}
{{\bf A}^X}(x,y) = (A_1^X(x,y),A_2^X(x,y))=\frac 12\, \left( -\frac{y-y_{0}}{r^2}, \frac{x-x_{0}}{r^2}\right)\,.
\end{equation}
The ${\Ab\Bb}X$-Hamiltonian is defined  by considering the Friedrichs
extension starting from $ C_0^\infty(\dot \Omega_{X})$
 and the associated differential operator is
\begin{equation}
-\Delta_{{\bf A}^X} := (D_x - A_1^X)^2 + (D_y-A_2^X)^2\,\mbox{with }D_x =-i\pa_x\mbox{ and }D_y=-i\pa_y.
\end{equation}
Let $ K_{X}$ be the antilinear operator 
$  K_{X} = e^{i \theta_{X}} \; \Gamma\,$,
with $  (x-x_0)+ i (y-y_0) = \sqrt{|x-x_0|^2+|y-y_0|^2}\, e^{i\theta_{X}}\,$,  $\theta_X$ such that $d\theta_X= 2 \Ab^X\,$, 
and where $ \Gamma$ is the complex conjugation operator  
$ \Gamma u = \bar u\,$.  A  function $ u$ is called  $ K_{X}$-real, if 
$ K_{X} u =u\,.$ The operator $ -\Delta_{\Ab^X}$ is preserving the 
$ K_{X}$-real functions and we can consider a
basis of $K_{X}$-real eigenfunctions. 
Hence we only  analyze the
restriction  of the ${\Ab\Bb}X$-Hamiltonian
to the $ K_{X}$-real space $ L^2_{K_{X}}$ where
$$
 L^2_{K_{X}}(\dot{\Omega}_{X})=\{u\in L^2(\dot{\Omega}_{X}) \;,\; K_{X}\,u =u\,\}\,.
$$
It was shown  in \cite{HHOO} and \cite{AFT} that the nodal set of such a $ K_X$ real eigenfunction has
the same structure as the nodal set of an eigenfunction of the
Laplacian except that an odd number of half-lines meet at $ X$. In particular, 
for a $K_X$-real groundstate (one pole), one can prove \cite{HHOO} that the nodal set consists of one line joining the pole and the boundary.

{\bf Extension to many poles}~\\
 We can extend this construction 
  in the case of a
 configuration with $ \ell$ distinct points $ X_1,\dots, X_\ell$ (putting a  flux $ \pi$ at each
 of these points). We just take as magnetic potential 
$$ 
\Ab^{\bf X} = \sum_{j=1}^\ell \Ab^{X_j}\,, \mbox{ where } \Xb=(X_1,\dots,X_\ell)\,.$$
We can also construct the antilinear
 operator $ K_{\bf X}$,  where $ \theta_X$ is replaced by a
 multivalued-function $ \phi_\Xb$ such that $ d\phi_X = 2 \Ab^{\Xb}$.  We can then  consider the 
 real subspace of the $ K_{\bf X}$-real
 functions in $ L^2_{K_{\Xb}}(\dot{\Omega}_{\Xb})$.
 It has been shown   that the $ K_{\Xb}$-real eigenfunctions have a regular nodal set
 (like the eigenfunctions of the Dirichlet Laplacian) with the
 exception that at each singular point $ X_j$ ($ j=1,\dots,\ell$) 
 an odd number of half-lines  meet. We recall  the following theorem which is the most interesting part of the magnetic characterization of the minimal partitions given in \cite{HHO3}:
\begin{theorem}\label{thchar}~\\
Let $ \Omega$ be simply connected. If $ \mathcal D$ is a $k$-minimal partition of $\Omega$, then, by choosing\footnote{We recall that $X(N)$ is defined after Definition \ref{AMS}.}  $ (X_1,\dots, X_\ell)= X^{odd}(N(\mathcal D))$, $ \mathcal D$ is the nodal partition 
of some $ k$-th $ K_{\Xb}$-real eigenfunction  of  the Aharonov-Bohm Laplacian  associated with $ \dot{\Omega}_{\bf X}$.
\end{theorem}

\subsection{Analysis of the critical sets in the large limit case} We first consider the case of one pole $X$. We look at a sequence of $K_X$-real eigenfunctions and follow the proof of Pleijel on the number of nodal domains. We observe that the part devoted to the lower bound works along the same lines and the way we shall meet $\mathfrak L_k(\Omega)$ is unchanged. When using  the Weyl formula, we observe that only a lower  bound of the counting function is used (see around \eqref{Wf}). If the distance of $X$ to the boundary is larger than $\epsilon$, we introduce a  disk $D(X,\epsilon)$ of radius 
  $\epsilon$ around $X$ ($\epsilon>0$) and consider the Dirichlet magnetic Laplacian in $\Omega \setminus \bar D (X,\epsilon)$.  For the $X$ at the distance less than $\epsilon$ of the boundary, we look at the magnetic Laplacian on $\Omega$ minus  a $(2\epsilon)$-tubular neighborhood of the boundary. In the two cases, we get an elliptic operator where the main term is the Laplacian $-\Delta$.  Hence we can combine the monotonicity of the Dirichlet problem with respect to the variation of the domain to the use of the standard Weyl formula (see \eqref{Wf}) to get (uniformly for $X$ in $\Omega$), an estimate
  for the counting function $N_{X} (\lambda)$ of $-\Delta_{\Ab^\Xb}$ in the following way:\\
 {\it  There exists a constant $C>0$ such that, for any $\epsilon >0$, as $\lambda \ar +\infty$,
 $$
  N_{X}(\lambda) \geq  \frac{1}{4 \pi} (1-C \epsilon) A(\Omega)\,\lambda +o(\lambda)\,.
  $$
  }
 Hence, for any $\epsilon >0$, any $X\in \Omega$, 
 $$
 \limsup_{n\ar +\infty} \mu(\phi_n^X)/n \leq (1+C \epsilon) \frac{4\pi }{ A(\Omega) \liminf_{k\ar +\infty} \frac{ \mathfrak L_k(\Omega)}{k}}\,.
 $$ 
 Taking the limit $\epsilon \ar 0$,  we get:
 \begin{equation}
 \limsup_{n \ar +\infty} \mu(\phi_n^X)/n \leq  \frac{4\pi }{ A(\Omega) \liminf_{k\ar +\infty}  \frac{ \mathfrak L_k(\Omega)}{k}}\,.
 \end{equation}
 Till now  $X$ was fixed. But everything being uniform with respect to $X$, we can also consider a sequence $\phi_n^{X_n}$ corresponding to the $n$-th eigenvalue of $- \Delta_{\Ab_{X_n}}$.\\
 Suppose that for a subsequence $k_j$, we have a $k_j$-minimal partition with only one pole $X_j$ in $\Omega$. Let $\phi_{k_j}^{X_j}$
  the corresponding eigenfunction. Hence, we are in a Courant sharp situation. The inequality above leads this time (possibly after extraction of a subsequence)
   to
   $$
   1 \leq  \frac{4\pi }{ A(\Omega) \liminf_{k\ar +\infty}  \frac{ \mathfrak L_k(\Omega)}{k}}\leq \nu_{Pl}\sim 0. 691\,.
 $$
 Hence a contradiction.\\
 We can play the same game with more than one pole and  get as consequence:
 \begin{proposition}
 If for  $k\in \mathbb N$, $\mathcal D_k$ denotes a minimal $k$-partition, then
  \begin{equation}
\lim_{k\ar + \infty}   \# X^{odd}(N(\mathcal D_k )) = +\infty\,.
\end{equation}
  \end{proposition}
  
 {\bf Proof.}\\
 Suppose indeed that this cardinality does not tend to $+\infty$. We can then extract a subsequence such that this cardinality  is finite.
 After  new  extractions of a subsequence, we can assume that this cardinality is fixed and that each critical point tends to a limiting point, which could be either at the boundary $\pa \Omega$ or in $\Omega$.  We apply Theorem~\ref{thchar} and consider the associated Aharonov-Bohm hamiltonians, whose poles are these odd critical points.
 We can then find a finite number of disks of radius $\epsilon$ centered at these limiting poles such that all the poles are contained in these balls for $k$  large enough. Then outside of these balls the potential $A^X$ and the derivatives are bounded by a uniform bound (depending on $\epsilon$) and the same construction works and leads to a contradiction.\\
 
 \begin{remark}
We recall that an upper bound for $ \# X (N(\mathcal D_k ))$  is given in \cite{HHO}  (case with no holes) by using Euler's formula:
  \begin{equation}\label{ubpc}
\# X^{odd}(N(\mathcal D_k )) \leq 2k - 4 \,. 
\end{equation}
On the other hand, the hexagonal conjecture for the asymptotic number of odd critical points of a $k$-minimal partition reads as follows:
\begin{equation}
\lim_{k\ar + \infty} \frac{\# X^{odd}(N(\mathcal D_k) ) }{k} = 2\,.
\end{equation}
Hence there are good reasons to believe 
that the upper bound \eqref{ubpc} 
is asymptotically optimal.
\end{remark}
 {\bf Acknowledgements.}\\
The discussions on this subject  started a few years ago with many other colleagues including M. Van den Berg, V. Bonnaillie-No\"el, G. Vial,  I. Polterovich, S. Steinerberger, .. and have continued during  various meetings (Oberwolfach, Vienna,  Rennes, Montreal, Banff, Loughborough, ...). Marie Helffer helped us for the drawing of the pictures and Corentin Lena transmitted to us enlightening pictures.\\

\end{document}